%% file: main.tex
\documentclass[conf]{new-aiaa}

\input{preamble/packages}

\input{preamble/commands}
\input{preamble/title}

\begin{document}

\maketitle

\input{sections/0-abstract}

\input{sections/1-notation}

\input{sections/2-intro}

\input{sections/3-ctcs-formulation}

\input{sections/4-solver}

\input{sections/5-pdg-formulation}

\input{sections/6-gpu}

\input{sections/7-numerical-results}

\input{sections/8-conclusion}



\section*{Acknowledgments}
This research was supported by NASA grant NNX17AH02A, AFOSR grant FA9550-20-1-0053, ONR grant N00014-20-1-2288, and Blue Origin, LLC; Government sponsorship is acknowledged.

\bibliography{references}

\end{document}

%% file: preamble/packages.tex
\usepackage{graphicx}
\usepackage{bm,amsmath}
\usepackage{xcolor}
\usepackage{amssymb}
\usepackage[version=4]{mhchem}
\usepackage{longtable,tabularx}
\usepackage{pgfplots}
\usepackage{textcomp}
\usepackage{subcaption}
\usepackage{placeins}

\pgfplotsset{compat=1.5}
\usetikzlibrary{positioning}
\usepackage{listings}
\definecolor{codegreen}{rgb}{0,0.6,0}
\definecolor{codegray}{rgb}{0.5,0.5,0.5}
\definecolor{codepurple}{rgb}{0.58,0,0.82}
\definecolor{backcolour}{rgb}{0.95,0.95,0.92}
\lstdefinestyle{mystyle}{
    backgroundcolor=\color{backcolour},   
    commentstyle=\color{codegreen},
    keywordstyle=\color{magenta},
    numberstyle=\tiny\color{codegray},
    stringstyle=\color{codepurple},
    basicstyle=\ttfamily\footnotesize,
    breakatwhitespace=false,         
    breaklines=true,                 
    captionpos=b,                    
    keepspaces=true,                 
    numbers=left,                    
    numbersep=5pt,                  
    showspaces=false,                
    showstringspaces=false,
    showtabs=false,                  
    tabsize=2
}

\usepackage{algorithm}
\usepackage{algpseudocode}

\usepackage{physics} 
\usepackage{mathtools} 
\usepackage{relsize} 
\usepackage{accents}

\hypersetup{
    colorlinks=true,
    linkcolor=blue!70,
    filecolor=PineGreen,
    urlcolor=red!50!gray!70,
    citecolor=red!50!gray!70,
    pdftitle={GPU SciTech 2024},
    pdfauthor={ACL}
}

\usepackage{tcolorbox}
\tcbuselibrary{skins}

\newtcolorbox{mybox}
{
  enhanced jigsaw,
  colframe=black,
  colback=white,
  drop shadow=black!50!white,
  boxrule=0.75pt
}

\newtcolorbox{mybox2}
{
  enhanced jigsaw,
  colframe=black,
  colback=white,
  drop shadow=black!50!white,
  boxrule=0.75pt,
  hbox
}

%% file: preamble/commands.tex
\newcommand{\pipg}{\textsc{pipg}}
\newcommand{\ecos}{\textsc{ecos}}
\newcommand{\osqp}{\textsc{osqp}}
\newcommand{\openblas}{\textsc{o}{\relsize{-1}pen}\textsc{blas}}
\newcommand{\pipgc}{\textsc{pipg}\textsubscript{{\relsize{-1}custom}}}
\newcommand{\cuda}{\textsc{cuda}}
\newcommand{\cublas}{{\relsize{-1}cu}\textsc{blas}}
\newcommand{\cusparse}{{\relsize{-1}cu}\textsc{sparse}}
\newcommand{\seco}{\textsc{s}{\relsize{-1}e}\textsc{co}}

\newcommand{\acikmese}{A\c{c}\i{}kme\c{s}e}
\newcommand{\behcet}{Beh\c{c}et}

\newcommand{\mc}[1]{\mathcal{#1}}
\newcommand{\mr}[1]{\mathrm{#1}}
\newcommand{\bR}[1]{\mathbb{R}^{#1}}
\newcommand{\derv}[1]{\overset{\circ}{#1}}

\let\norm\undefined
\newcommand{\norm}[1]{\big\|#1\big\|}

%% file: preamble/title.tex
\title{\Large Fast Monte Carlo Analysis for 6-DoF Powered-Descent Guidance\\via GPU-Accelerated Sequential Convex Programming}

\author{Govind M.\ Chari$^*$\!, Abhinav G.\ Kamath$^*$\!, Purnanand Elango\footnote{Ph.D.\ Student, William E.\ Boeing Department of Aeronautics \& Astronautics; \texttt{\{gchari,\,agkamath,\,pelango\}@uw.edu}}\!, and \behcet~\acikmese\footnote{Professor, William E.\ Boeing Department of Aeronautics \& Astronautics; \texttt{behcet@uw.edu}}}
\affil{University of Washington, Seattle, WA 98195, USA}

%% file: sections/0-abstract.tex
\begin{abstract}
We introduce a GPU-accelerated Monte Carlo framework for nonconvex, free-final-time trajectory optimization problems. This framework makes use of the prox-linear method, which belongs to the larger family of sequential convex programming (SCP) algorithms, in conjunction with a constraint reformulation that guarantees inter-sample constraint satisfaction. Key features of this framework are: (1) continuous-time constraint satisfaction; (2) a matrix-inverse-free solution method; (3) the use of the proportional-integral projected gradient (PIPG) method, a first-order convex optimization solver, customized to the convex subproblem at hand; and, (4) an end-to-end, library-free implementation of the algorithm. We demonstrate this GPU-based framework on the 6-DoF powered-descent guidance problem, and show that it is faster than an equivalent serial CPU implementation for Monte Carlo simulations with over 1000 runs. To the best of our knowledge, this is the first GPU-based implementation of a general-purpose nonconvex trajectory optimization solver.
\end{abstract}

%% file: sections/1-notation.tex


\section*{Notation}

\begin{tabular}{ll}
$\bR{n}\!=\!\bR{n\times 1}$ & Vectors are matrices with one column\\
$1_n,0_n$   & Vector of ones and zeros, respectively, in $\bR{n}$ (subscript inferred whenever omitted)\\
$0_{n\times m}$ & Matrix of zeros in $\bR{n\times m}$\\
$I_n$   & Identity matrix in $\bR{n\times n}$\\
$V_{[1:n]}$ & Collection of vectors or matrices $V_k$, for $k=1,\ldots,n$\\
$|v|_+$ & $\square \mapsto \max\{0,\square\}$ applied element-wise for vector $v$\\
$v^2$   & $\square \mapsto \square^2$ applied element-wise for vector $v$\\
$(u,v)$ & Concatenation of vectors $u\in\bR{n}$ and $v\in\bR{m}$ to form a vector in $\bR{n+m}$\\
$[A\,B]$ & Concatenation of matrices $A$ and $B$ with same number of rows\\
$\bigg[\!\!\!\begin{array}{c}A\\[-0.1cm]B\end{array}\!\!\bigg]$ & Concatenation of matrices $A$ and $B$ with same number of columns\\
\end{tabular}

%% file: sections/2-intro.tex
\section{Introduction}\label{sec:intro}

Monte Carlo analysis is critical for certifying the safety and performance of guidance, navigation, and control (GNC) systems \cite{malyuta2019discretization}. For a 6-DoF powered-descent guidance algorithm, Monte Carlo analysis would involve generating hundreds or even thousands of dispersed initial conditions around the nominal initial condition, computing trajectories with these initial conditions, and assessing statistics such as solve-times, percentage of converged trajectories, and propellant consumption. However, generating thousands of trajectories serially can be time consuming, and can lead to increased development times. Since Monte Carlo simulation is \textit{embarrassingly parallel} \cite{herlihy2020art}, we propose performing Monte Carlo analyses on a graphics processing unit (GPU) to exploit its parallel computing capability.

Due to their parallel computing capability and high memory bandwidth, GPUs have become increasingly popular in fields that require solving large-scale problems (in which parallelism can be exploited). Most notably, GPUs have become a staple in machine learning for training deep learning models. For example, OpenAI has a supercomputer with 10,000 GPUs to train its large language models \cite{openai-gpu}, and Tesla, reportedly, also has a 10,000-GPU cluster to perform training for its machine learning-based full self-driving technology \cite{Shilov_2023}. Additionally, many modern machine learning tools such as JAX, TensorFlow, and PyTorch, provide programmers with user-friendly, high-level interfaces to leverage the compute power of GPUs \cite{tensorflow2015-whitepaper, pytorch, jax2018github}.

GPUs have also been used outside of machine learning for a variety of tasks. For example, \cite{Ilg2011} uses a GPU for running projectile simulations, i.e., propagating the equations of motion for a projectile, in parallel.
In \cite{Rastgar2023}, the authors develop a framework to perform trajectory optimization on GPUs; however, they consider a restrictive trajectory optimization problem template that does not model dynamics or generic nonconvex path constraints. The only constraints considered are maximum speed, maximum acceleration, and keepout zone constraints. As a result, the generated trajectories are not guaranteed to obey the equations of motion of the vehicle. 

In this paper, we demonstrate a GPU-based parallelized implementation for Monte Carlo analysis of trajectory optimization problems with general nonlinear dynamics and nonconvex path constraints. The solutions obtained are guaranteed to be feasible with respect to the original nonlinear dynamics, and all the imposed constraints are guaranteed to be satisfied at all points in time, i.e., in continuous-time (subject to successful convergence and up to numerical integration and optimization tolerances). This final guarantee on continuous-time constraint satisfaction is a result of the problem formulation proposed in \cite{ctcs2024}. To demonstrate this framework, we consider the 6-DoF powered-descent guidance algorithm \cite{miki2020successive}.

The propellant-optimal powered-descent guidance (PDG) problem is concerned with finding a feasible sequence of control inputs that allows a vehicle to decelerate itself from some initial position and velocity, to perform a soft landing at some specified target on the ground, while satisfying a variety of state and control constraints. This problem, in general, is of great importance to companies like SpaceX that land their boosters propulsively, allowing them to make launches cheaper \cite{blackmore2016autonomous}, and for NASA JPL, to land their rovers near important scientific targets on Mars \cite{martin2015flight}, for example. Solving the 3-DoF PDG problem, where only translational dynamics are considered, is one approach to solving this problem. However, since this method does not explicity consider attitude dynamics, attitude constraints, such as line-of-sight, cannot be enforced directly. For applications such as NASA's Safe and Precise Landing — Integrated Capabilities Evolution (SPLICE) project, the terrain relative navigation (TRN) cameras and hazard detection LiDAR (HDL) require constraining the attitude of the rocket, necessitating the use of 6-DoF PDG algorithms \cite{carson2019splice}.

General trajectory optimization problems, such as the 6-DoF PDG problem, can be posed as nonconvex optimization problems in continuous-time \cite{malyuta2021advances}. Given the nonconvexity of these problems, it is typically prohibitively expensive, if not outright impossible, to obtain globally optimal trajectories. As a result, we must turn to solution procedures, such as sequential convex programming (SCP) \cite{SCPToolboxCSM2022}, that yield locally optimal solutions. At a high level, SCP involves (1) convexification (often via linearization) of all nonconvex constraints about some reference trajectory to yield an infinite-dimensional convex optimization problem; (2) discretization, to convert the problem to a numerically tractable finite-dimensional convex optimization problem; and finally (3) finding a solution to this finite-dimensional convex subproblem to update the reference trajectory. This process is applied iteratively until convergence. As the name suggests, this process is inherently sequential in nature. However, we can leverage parallelism to simultaneously solve for many trajectories in Monte Carlo analysis.

There are a number of factors that make developing a GPU-based parallelized implementation of SCP challenging. SCP requires the use of many linear algebra operations and the solution to a sequence of convex optimization subproblems. For a CPU-based implementation, we can use linear algebra packages such as {\openblas} \cite{openblas} and open-source convex optimization solvers such as {\osqp} and {\ecos} \cite{osqp, domahidi2013ecos}. However, {\openblas} binaries are compiled for CPU architectures and NVIDIA's linear algebra libraries, {\cublas} and {\cusparse}, are written such that functions from these libraries must be called by the CPU, even though the computations are performed on the GPU. This is undesirable for our implementation of SCP, since we require all computations for a single trajectory to be performed by a single GPU thread. Similarly, although a GPU-accelerated implementation of {\osqp} exists, the solver itself must be called by the CPU \cite{osqp-gpu}. Although off-the-shelf convex optimization solvers could be modified to work on the GPU, custom parser frameworks are typically required—to put the convex subproblem in the canonical form prior to being passed to the solver—as well \cite{lofberg2004yalmip, diamond2016cvxpy}. Thus, for our GPU-based implementation of SCP, we implement Level 1, 2, and 3 BLAS operations for discretization and other such operations, and the first-order proportional-integral projected gradient ({\pipg}) solver \cite{yu2022extrapolated, kamath2023seco}, customized to the trajectory optimization problem template \cite{elango2022customized, kamath2023customized}, allowing us to bypass parsing entirely. Further, as in the sequential conic optimization ({\seco}) framework described in \cite{kamath2023customized}, this framework completely eliminates the need for sparse linear algebra operations, is entirely matrix-factorization/inversion-free, and is implemented using a library-free codebase.

%% file: sections/3-ctcs-formulation.tex
\section{Problem Formulation}\label{sec:ctcs}
We consider a continuous-time dynamical system of the form:
\begin{align}
    \frac{\mr{d}\xi(t)}{\mr{d}t} = \dot{\xi}(t) = F(\xi(t),\zeta(t))\label{sys-dyn}
\end{align}
for $t\in[0,t_{\mr{f}}]$, with state $\xi\in\bR{n_\xi}$ and control input $\zeta\in\bR{n_\zeta}$, which are subject to path constraints:
\begin{subequations}
 \begin{align}
    & g(\xi(t),\zeta(t)) \le 0_{n_g}\\
    & h(\xi(t),\zeta(t)) = 0_{n_h}
\end{align}\label{path-cnstr}%
\end{subequations}
where $g$ and $h$ are vector-valued functions. The operators ``$\le$'' and ``$=$'' are interpreted element-wise. Boundary conditions at the initial and final times could also be specified.

The goal of trajectory optimization is to generate a solution $\xi(t),\zeta(t)$ for \eqref{sys-dyn} that satisfies the path constraints, boundary conditions, and optimizes a linear function of the terminal state $\xi(t_{\mathrm{f}})$, i.e., solve the following Mayer form optimal control problem \cite[Sec. 3.3.2]{liberzon2011calculus}:
\begin{subequations}
\begin{align}
\underset{t_{\mathrm{f}},\,\xi,\,\zeta}{\mathrm{minimize}}~~&~\xi(t_{\mathrm{f}})^\top e_{\mathrm{cost}} & & \\
\mathrm{subject~to}~&~\dot{\xi}(t) = F(\xi(t),\zeta(t)) & & t\in[0,t_{\mathrm{f}}] \\
 &~g(\xi(t),\zeta(t)) \le 0 & & t\in[0,t_{\mathrm{f}}] \\
 &~h(\xi(t),\zeta(t)) = 0 & & t\in[0,t_{\mathrm{f}}] \\
 &~E_{\mathrm{i}}\xi(0) = z_{\mathrm{i}},~E_{\mathrm{f}}\xi(t_{\mathrm{f}}) = z_{\mathrm{f}}
\end{align}
\end{subequations}
where $E_{\mathrm{i}}\in\bR{n_\mathrm{i}\times n_{\xi}}$ and $E_{\mathrm{f}}\in\bR{n_\mathrm{f}\times n_{\xi}}$ select components of the initial and final state, respectively.
\subsection{Time-Dilation and Constraint Reformulation}
When the final time, $t_{\mathrm{f}}$, is an unknown, we adopt a well-known technique called 
time-dilation \cite{ctcs2024,kamath2023seco,szmuk2020successive} to convert a free-final-time problem to an equivalent fixed-final-time problem: 
\begin{align}
     \frac{\mr{d}\xi(t(\tau))}{\mr{d}\tau} = \derv{\xi}(\tau) = \frac{\mr{d}\xi(\tau)}{\mr{d}t}s(\tau) = s(\tau)F(\xi(\tau),\zeta(\tau))\label{sys-dyn-dil}
\end{align}
for $\tau\in[0,1]$, where $t(\tau)$ is a strictly increasing function with domain $[0,1]$. Note that $\xi$ and $\zeta$ are re-defined to be functions of $\tau\in[0,1]$ instead of $t$. The dilation factor: 
$$
    s(\tau) = \frac{\mr{d}t(\tau)}{\mr{d}\tau}
$$ 
is treated as an additional control input and is required to be positive. The control input is then re-defined to be $u = (\zeta,s)$, and the final time is given by:
$$
    t_{\mr{f}} = \int_0^1s(\tau)\mr{d}\tau
$$

Next, to address the long-standing issue of inter-sample constraint violation in direct methods \cite{betts2010practical}, we augment the path constraints to the time-dilated system dynamics \eqref{sys-dyn-dil} as follows:
\begin{align}
    & \derv{x}(\tau) = \begin{bmatrix}
                           \derv{\xi}(\tau) \\
                           \derv{y}(\tau) 
                       \end{bmatrix} = s(\tau)\begin{bmatrix}
                                           F(\xi(\tau),\zeta(\tau)) \\[0.15cm]
                                           1^\top |g(\xi(\tau),\zeta(\tau))|^2_+ + 1^\top h(\xi(\tau),\zeta(\tau))^2
                                       \end{bmatrix} = f(x(\tau),u(\tau))\label{sys-dyn-dil-aug}
\end{align}
for $\tau\in[0,1]$, where $x=(\xi,y)$ is the re-defined state. Periodic boundary condition on the constraint violation integrator:
\begin{align}
    y(0) = y(1) \label{cnstr-integrator-bc}
\end{align}
ensure that the path constraints are satisfied for all $\tau\in[0,1]$. We refer the reader to \cite{ctcs2024,sargent1977development} for further details.
\subsection{Discretization}
Discretization is a crucial step in direct numerical optimal control for computational tractability, i.e., for obtaining an approximate finite-dimensional optimization problem from the original infinite-dimensional optimal control problem. The first step is to generate a grid with $N$ nodes within the interval $[0,1]$:
\begin{align}
    0 = \tau_1 < \ldots < \tau_N = 1\label{disc-grid}
\end{align}
The time duration over sub-interval $[\tau_k,\tau_{k+1}]$ is given by:
\begin{align*}
\Delta t_k = \int_{\tau_k}^{\tau_{k+1}} s(\tau)\mr{d}\tau
\end{align*}
for $k=1,\ldots,N-1$, and the following constraints are specified on the dilation factor and time duration:
\begin{align*}
    s_{\min} \le s(\tau) \le{} & s_{\max} \\
    t_{\mr{f}} = \sum_{k=1}^{N-1}\Delta t_k \le{} & t_{\mr{f},\max}\\    
    \Delta t _{\min} \le \Delta t_k \le{} & \Delta t_{\max}
\end{align*}

Next, we parameterize the control input with a first-order-hold (FOH):
$$
    u(\tau) = \left(\frac{\tau_{k+1}-\tau}{\tau_{k+1}-\tau_k}\right)u_k + \left(\frac{\tau-\tau_{k}}{\tau_{k+1}-\tau_k}\right)u_{k+1}
$$
for $\tau\in[\tau_{k},\tau_{k+1}]$ and $k=1,\ldots,N-1$. 
In practice, the combination of time-dilation and FOH parameterization for control input is sufficient for accurately modeling a large class of optimal control problems. See \cite{malyuta2019discretization} for a detailed discussion on the relative merits of popular discretization and parameterization techniques in numerical optimal control.

We treat the node-point values of the state and control input, $x_k$ and $u_k$, respectively, as decision variables, and re-write the system dynamics \eqref{sys-dyn-dil-aug} equivalently as:
\begin{align}
    x_{k+1} = x_{k} + \int_{\tau_k}^{\tau_{k+1}}f(x(\tau),u(\tau))\mr{d}\tau\label{disc-dyn}
\end{align}
for $k=1,\ldots,N-1$, where $x(\tau)$ is the solution to \eqref{sys-dyn-dil-aug} over $[\tau_k,\tau_{k+1}]$ with initial condition $x_k$, and $u(\tau)$ is FOH-parameterized using $u_k$ and $u_{k+1}$. Then, \eqref{disc-dyn} can be compactly expressed through function $f_k:\bR{n_x}\times\bR{n_x}\times\bR{n_u}\times\bR{n_u}\to\bR{n_x}$ as follows:
\begin{align}
    x_{k+1} = f_k(x_{k},u_{k},u_{k+1})
\end{align}

The boundary condition \eqref{cnstr-integrator-bc} is relaxed to: 
$$
    y_{k+1} - y_{k} \le \epsilon
$$
for $k=1,\ldots,N-1$ and a sufficiently small $\epsilon>0$, so that all feasible solutions do not violate LICQ \cite{ctcs2024}. Note that $y_k$ is the value of the state $y$ in \eqref{sys-dyn-dil-aug} at node $\tau_k$.

The discretized optimal control problem with parameterized control input is given by:
\begin{subequations}
\begin{align}
\underset{x_k,u_k}{\mathrm{minimize}}~~&~x_N^\top \tilde{e}_{\mathrm{cost}} & & \\
\mathrm{subject~to}~&~x_{k+1} = f_k(x_{k},u_k,u_{k+1}) & & k=1,\ldots,N-1\label{disc-ocp:dyn}\\\
 &~E_y (x_{k+1} - x_k) \le \epsilon & & k=1,\ldots,N-1\\
 &~u_{\min} \le u_k \le u_{\max} & & k=1,\ldots,N\\
 &~\tilde{E}_{\mathrm{i}}x_1 = (z_{\mathrm{i}},0),~\tilde{E}_{\mathrm{f}}x_N = z_{\mathrm{f}}
\end{align}\label{disc-ocp}%
\end{subequations}
where
\begin{align}
    & \tilde{e}_{\mathrm{cost}} = (e_{\mathrm{cost}},0)\\
    & E_y = [0_{1\times n_\xi}~~1]\\
    & \tilde{E}_{\mathrm{i}} = \begin{bmatrix}
                                    E_{\mathrm{i}} & 0_{n_\xi}\\
                                    0_{1\times n_\xi} & 1
                               \end{bmatrix}\\
    & \tilde{E}_{\mathrm{f}} = \big[E_{\mathrm{f}}~0_{n_\xi}\big]
\end{align}

\subsection{Sequential Convex Programming}

We adopt the prox-linear method \cite{drusvyatskiy2018error}, an SCP algorithm, to solve:
\begin{subequations}
\begin{align}
    \underset{x_k,u_k}{\mathrm{minimize}}~~&~x_N^\top\tilde{e}_{\mathrm{cost}} + w_{\mr{ep}}\sum_{k=1}^{N-1}\|x_{k+1}-f_k(x_k,u_k,u_{k+1})\|_1 & & \\
   \mathrm{subject~to}~&~E_y(x_{k+1} - x_k) \le \epsilon & & \label{disc-ocp-expen:relax}\\
                       &~u_{\min} \le u_k \le u_{\max} & & \label{disc-ocp-expen:ctrl}\\
                       &~\tilde{E}_{\mathrm{i}}x_1 = (z_{\mathrm{i}},0),~\tilde{E}_{\mathrm{f}}x_N = z_{\mathrm{f}}\label{disc-ocp-expen:bc}
\end{align}\label{disc-ocp-expen}%
\end{subequations}
where the (only) nonconvex constraint of \eqref{disc-ocp} is exactly penalized with $\ell_1$-norm. For a large enough, finite $w_{\mr{ep}}$, we can show that a solution of \eqref{disc-ocp-expen} that is feasible with respect to \eqref{disc-ocp:dyn} is also a solution of \eqref{disc-ocp}. See \cite[Chap. 17]{nocedal2006numerical} for a detailed discussion.

The penalized problem \eqref{disc-ocp-expen} can be compactly stated as:
\begin{align}
    \underset{Z\in\mathcal{Z}}{\mathrm{minimize}}~~L(Z) + w_{\mr{ep}}\|G(Z)\|_1\label{disc-ocp-expen-compact}
\end{align}
where $Z = (x_1,\ldots,x_N,u_1,\ldots,u_N)$ and $\mathcal{Z}$ is a convex set corresponding to constraints \eqref{disc-ocp-expen:relax}, \eqref{disc-ocp-expen:ctrl}, and \eqref{disc-ocp-expen:bc}.

The prox-linear method determines a stationary point of \eqref{disc-ocp-expen-compact} by solving a sequence of convex subproblems. At iteration $j+1$, it solves the following convex problem:
\begin{align}
    \underset{Z\in\mathcal{Z}}{\mathrm{minimize}}~~L(Z) + w_{\mr{ep}}\|G(Z^{j})+\nabla G(Z^{j})(Z-Z^{j})\|_1 + \frac{w_{\mr{prox}}}{2}\|Z-Z^j\|_2^2\label{cvx-subproblem}
\end{align}
where $Z^j$ is the solution from iteration $j$. The sequence $Z^j$ can be shown to converge for an appropriate choice of the proximal term weight $w_{\mr{prox}}$. Further, if the converged point is feasible for \eqref{disc-ocp}, then it is also a KKT point.

There are two important steps in the construction of each convex subproblem--(1) computation of $\nabla G$: the gradient of $f_k$ with respect to the previous iterate, and (2) prescaling and preconditioning.

\subsubsection{Gradient of Discretized Dynamics}

At each iteration of SCP, the discretized dynamics constraint \eqref{disc-ocp:dyn} is linearized with respect to the previous iterate, denoted by $\bar{x}_k,\,\bar{u}_k$, for $k=1,\ldots,N$. 

For each $k=1,\ldots,N-1$, let $\bar{x}^k(\tau)$ denote the solution to \eqref{sys-dyn-dil-aug} over $[\tau_k,\tau_{k+1}]$ generated with initial condition $\bar{x}_k$ and control input $\bar{u}(\tau)$, which is parameterized using $\bar{u}_k$, for $k=1,\ldots,N$. Note that, in general, $\bar{x}^k(\tau_{k+1})\ne \bar{x}_{k+1}$ before SCP converges \cite{kamath2023seco}.

The Jacobians of $f$ in \eqref{sys-dyn-dil-aug}, evaluated with respect to $\bar{x}^k(\tau),\,\bar{u}(\tau)$, are denoted by:
\begin{align*}
    A^k(\tau) ={} & \frac{\partial f(\bar{x}^k(\tau),\bar{u}(\tau))}{\partial x}\\
    B^k(\tau) ={} & \frac{\partial f(\bar{x}^k(\tau),\bar{u}(\tau))}{\partial u}
\end{align*}
for $\tau\in[\tau_k,\tau_{k+1}]$.

The discretized dynamics constraint \eqref{disc-dyn} is linearized as follows:
\begin{align}
    x_{k+1} = A_kx_{k} + B^-_ku_k + B^+_ku_{k+1} + w_k     
\end{align}
for $k=1,\ldots,N-1$, where $A_k$, $B^-_k$, $B^+_k$, $w_k$ result from the solution to the following initial value problem over $\tau\in[\tau_k,\tau_{k+1}]$:
\begin{align*}
    & \derv{\Phi}_x(\tau,\tau_k) = A^k(\tau)\Phi_x(\tau,\tau_k)\\ 
    & \derv{\Phi}{}^-_u(\tau,\tau_k) = A^k(\tau)\Phi{}^-_u(\tau,\tau_k) + B^k(\tau)\left(\frac{\tau_{k+1}-\tau}{\tau_{k+1}-\tau_k}\right)\\
    & \derv{\Phi}{}^+_u(\tau,\tau_k) = A^k(\tau)\Phi{}^+_u(\tau,\tau_k) + B^k(\tau)\left(\frac{\tau-\tau_k}{\tau_{k+1}-\tau_k}\right)\\
    & \Phi_x(\tau_k,\tau_k) = I_{n_x}\\
    & \Phi{}^-_u(\tau_k,\tau_k) = 0_{n_x\times n_u}\\
    & \Phi{}^+_u(\tau_k,\tau_k) = 0_{n_x\times n_u}\\[0.5cm]
    & A_k = \Phi_x(\tau_{k+1},\tau_k)\\
    & B_k^- = \Phi{}^-_u(\tau_{k+1},\tau_k)\\
    & B_k^+ = \Phi{}^+_u(\tau_{k+1},\tau_k)\\
    & w_k = \bar{x}^k(\tau_{k+1}) - A_k\bar{x}_k - B_k^-\bar{u}_{k} - B_{k}^+\bar{u}_{k+1} 
\end{align*}

\subsubsection{Convex Subproblem with Prescaling and Preconditioning}

The final step in the construction of the convex subproblem involves scaling the decision variables and preconditioning the constraint parameters. In particular, we treat $\hat{x}_k$ and $\hat{u}_k$ as the decision variables, which are linearly related to $x_k$ and $u_k$ as follows: 
\begin{subequations}
\begin{align}
& x_k = P_x \hat{x}_k + \bar{x}_k\\
& u_k = P_u \hat{u}_k + \bar{u}_k
\end{align}
\end{subequations}
where $P_x$ and $P_u$ are selected to ensure that $\hat{x}_k$ and $\hat{u}_k$, for $k=1,\ldots,N$, are roughly on the same order of magnitude, which greatly aids in convergence. Next, we define:
\begin{subequations}
\begin{align}
        & \hat{A}^-_k = P_x^{-1}A_kP_x\\
        & \hat{A}^+_k = -I_{n_{x}}\\
        & \hat{B}^-_k = P_x^{-1}B^-_kP_u\\
        & \hat{B}^+_k = P_x^{-1}B^+_kP_u\\
        & \hat{w}_k = P_x^{-1}(\bar{x}^{k}(\tau_{k+1})-\bar{x}_{k+1})
\end{align}
\end{subequations}
Finally, the convex subproblem \eqref{cvx-subproblem}, which is a quadratic program (QP), can be stated as:
\begin{subequations}
\begin{align}
    \underset{\hat{x}_k,\hat{u}_k}{\mr{minimize}}~~&~w_{\mr{cost}}\hat{x}_N^\top e_{\mathrm{cost}}+ \frac{w_{\mr{prox}}}{2}\sum_{k=1}^N \|\hat{x}_k\|_2^2 + \|\hat{u}_k\|_2^2 + w_{\mr{ep}}\sum_{k=1}^{N-1}1_{n_x}^\top(\mu^+_k+\mu_k^-) & &  \label{cvx-subproblem-custom:obj}\\
    \mr{subject~to}~&~\hat{A}^-_k\hat{x}_k +\hat{A}_k^+\hat{x}_{k+1} + \hat{B}^-_k\hat{u}_k + \hat{B}^+_k\hat{u}_{k+1} + \mu_k^+ - \mu_k^- + \hat{w}_k = 0 & & \hspace{-1cm}1 \le k \le N-1\label{cvx-subproblem-custom:dyn}\\
     &~E_{y}(\hat{x}_{k+1}-\hat{x}_{k}) \le \hat{\varepsilon}_k & & \hspace{-1cm}1 \le k \le N-1\label{cvx-subproblem-custom:relax}\\
     &~\mu_k^+\ge 0,~\mu^-_k\ge 0 & &  \hspace{-1cm}1 \le k \le N-1 \\
     &~\hat{u}_{k,\min} \le \hat{u}_k \le \hat{u}_{k,\max} & &  \hspace{-1cm}1 \le k \le N\\
     &~\tilde{E}_{\mr{i}}\hat{x}_1 = \hat{z}_{\mr{i}},~\tilde{E}_{\mr{f}}\hat{x}_N = \hat{z}_{\mr{f}}
\end{align}\label{cvx-subproblem-custom}%
\end{subequations}
where
\begin{align*}
        & \hat{\varepsilon}_k = \epsilon E_yP_x^{-1}E_y^\top - E_yP_x^{-1}(\bar{x}_{k+1}-\bar{x}_k)\\
        & \hat{u}_{k,\min\!/\!\max} = P_u^{-1}(u_{\min\!/\!\max}-\bar{u}_k)\\
        & \hat{z}_{\mr{i}} = \tilde{E}_{\mr{i}}P_x^{-1}\tilde{E}_{\mr{i}}^\top((z_{\mr{i}},0) - \tilde{E}_{\mr{i}}\bar{x}_1)\\
        & \hat{z}_{\mr{f}} = \tilde{E}_{\mr{f}}P_x^{-1}\tilde{E}_{\mr{f}}^\top(z_{\mr{f}} - \tilde{E}_{\mr{f}}\bar{x}_N)
\end{align*}

%% file: sections/4-solver.tex
\section{Convex Optimization Solver}
The convex subproblem \eqref{cvx-subproblem-custom} is a quadratic program (QP), which can be represented in the canonical form as:
\begin{subequations}
\begin{align}
    \mr{minimize}~~&~\frac{1}{2}z^\top\hat{P} z + \hat{p}^\top z\\
    \mr{subject~to}~&~\hat{G}z = \hat{g}\\
    &~\hat{H}z\le \hat{h} 
\end{align}\label{canonical-QP}
\end{subequations}
where 
\begin{subequations}
\begin{align}
    & z = \big(\hat{x}_1,\ldots,\hat{x}_N,\,\hat{u}_1,\ldots,\hat{u}_N,\,\mu^-_1,\ldots,\mu^-_{N-1},\,\mu^+_1,\ldots,\mu^+_{N-1}\big)\in\bR{(n_x+n_u)N+2n_x(N-1)}\\
    & \hat{P} = \mr{blkdiag}(w_{\mr{tr}}I_{(n_x+n_u)N},\,0_{2n_x(N-1)\times 2n_x(N-1)})\\
    & \hat{p} = w_{\mr{cost}}\big(0_{n_x(N-1)},\,e_x,\,0_{n_uN+2n_x(N-1)}\big)\\
    & ~~~~~+w_{\mr{vc}}\big(0_{(n_x+n_u)N},\,1_{2n_x(N-1)}\big)\\
    & \hat{G} = \begin{bmatrix}\hat{G}_x & \hat{G}_u & \hat{G}_\mu\end{bmatrix}\\
    & \hat{G}_x = \left[\hat{A}^{-}~0_{n_x(N-1)\times n_x}\right]+\left[0_{n_x(N-1)\times n_x}~\hat{A}^{+}\right]\\
    & \hat{G}_u = \left[\hat{B}^-~0_{n_x(N-1)\times n_u}\right]+\left[0_{n_x(N-1)\times n_u}~\hat{B}^+\right]\\
    & \hat{G}_\mu = \left[I_{n_x(N-1)}~-I_{n_x(N-1)}\right]\\
    & \hat{A}^- = \mr{blkdiag}(\hat{A}^-_1,\ldots,\hat{A}^-_{N-1})\\
    & \hat{A}^+ = \mr{blkdiag}(\hat{A}^+_1,\ldots,\hat{A}^+_{N-1})\\
    & \hat{B}^- = \mr{blkdiag}(\hat{B}^-_1,\ldots,\hat{B}^-_{N-1})\\
    & \hat{B}^+ = \mr{blkdiag}(\hat{B}^+_1,\ldots,\hat{B}^+_{N-1})\\
    & \hat{g} = -\big(\hat{w}_1,\ldots,\hat{w}_{N-1}\big)\\
    & \hat{H}_y = -\left[I_{N-1}\otimes E_y~0_{n_y(N-1)\times n_x}\right] + \left[0_{n_y(N-1)\times n_x}~I_{N-1}\otimes E_y\right]\\
    & \hat{H} = \begin{bmatrix} \hat{H}_y~~0_{n_y(N-1)\times n_uN+2n_x(N-1)} \end{bmatrix}\\
    & \hat{h} = \big( \hat{u}_{1,\max},\ldots,\hat{u}_{N,\max},-\hat{u}_{1,\min},\ldots,-\hat{u}_{N,\min},\,0_{2n_x(N-1)},\,\hat{\varepsilon}_1,\ldots,\hat{\varepsilon}_{N-1}\big)
\end{align}\label{parsed-data}%
\end{subequations}
We adopt the proportional-integral projected gradient ({\pipg}) algorithm \cite{yu2022extrapolated} with extrapolation (shown in Algorithm \ref{alg:pipg}) to solve \eqref{canonical-QP}. The input to the algorithm consists of maximum iterations $j_{\max}$, the extrapolation parameter $\rho$, the step-size $\alpha$, which depends on the singular values of $\hat{P}$ and $\hat{H}$, and the step-size ratio $\omega$, which is chosen heuristically to balance the convergence rates of the primal and dual iterates.
\begin{algorithm}[!htpb]
\caption{{\pipg}}\label{alg:pipg}
\begin{flushleft}
\textbf{Inputs:}~$j_{\max},\,\rho,\,\alpha,\,\omega$
\end{flushleft}
\begin{algorithmic}[1]
\Statex \hspace{-0.6cm}\textbf{Initialize:} $\zeta^1,\,\eta^1,\,\chi^1$
\State $\beta\gets\omega\alpha$
\For{$j = 1,\ldots,j_{\max}$}
\State $z^{j+1} \gets \zeta^{j} - \alpha(\hat{P}\zeta^j+\hat{p}+\hat{G}^\top\eta^j + \hat{H}^\top\chi^j)$
\State {$\tilde{E}_{\mr{i}}\hat{x}_{1}^{j+1} \gets \hat{z}_{\mr{i}}$}
\State {$\tilde{E}_{\mr{f}}\hat{x}_{K}^{j+1} \gets \hat{z}_{\mr{f}}$}
\State {$\hat{u}_k^{j+1} \gets \max\{\hat{u}_{k,\min},\min\{\hat{u}_{k,\max},\hat{u}_k^{j+1}\}\}~\qquad k=1,\ldots,N$}
\State {$(\mu^{+}_k)^{j+1}\gets\max\{0,(\mu^{+}_k)^{j+1}\}~\quad\,\qquad\qquad\qquad k=1,\ldots,N-1$}
\State {$(\mu^{-}_k)^{j+1}\gets\max\{0,(\mu^{-}_k)^{j+1}\}~\quad\,\qquad\qquad\qquad k=1,\ldots,N-1$}
\State $w^{j+1} \gets \eta^j + \beta(\hat{G}(2z^{j+1}-\zeta^j)-\hat{g})$
\State $v^{j+1} \gets \max\{0,\chi^j + \beta(\hat{H}(2z^{j+1}-\zeta^j)-\hat{h})\}$
\State $\zeta^{j+1} \gets (1-\rho)\zeta^j + \rho z^{j+1}$
\State $\eta^{j+1} \gets (1-\rho)\eta^j + \rho w^{j+1}$
\State $\chi^{j+1} \gets (1-\rho)\chi^j + \rho v^{j+1}$
\EndFor
\end{algorithmic}
\begin{flushleft}
\textbf{Return:}~$z^{j},\,w^j,\,v^j$
\end{flushleft}
\end{algorithm}

Algorithm \ref{alg:pipg} can be customized to Algorithm \ref{alg:pipg_custom} (referred to as {\pipgc}), which exploits the sparsity structure of the parameters \eqref{parsed-data} of \eqref{canonical-QP}. The recursive calls to {\pipgc} within SCP are warm-started with the primal-dual solution from the previous subproblem to accelerate its convergence.

The customized power iteration method shown in Algorithm \ref{alg:power_custom}, called before {\pipgc}, estimates the square of the maximum singular value of $\hat{H}$, denoted by $\sigma$. The maximum eigenvalue of $\hat{P}$ is trivially known, given its construction as a diagonal matrix. Furthermore, {\pipgc} employs a customized termination criterion, as shown in Algorithm \ref{alg:stopping_custom}.

The primary design consideration in Algorithms \ref{alg:power_custom}, \ref{alg:pipg_custom}, and \ref{alg:stopping_custom} is to eliminate matrix factorizations, inversions, and sparse linear algebra; only simple operations on small-sized dense vectors and matrices are used. See \cite{yu2022extrapolated,elango2022customized,kamath2023customized} for further details on custom, real-time implementations of {\pipg}.

\begin{algorithm}[H]
\small
\caption{Customized power iteration method}
\label{alg:power_custom}
    \vspace{0.25em}
    \begin{flushleft}
        \textbf{Inputs:} $\hat{A}^{-}_{[1:N-1]}$, $\hat{A}^{+}_{[1:N-1]}$, $\hat{B}^{-}_{[1:N-1]}$, $\hat{B}^{+}_{[1:N-1]}$, $\hat{x}_{[1:N]}$, $\hat{u}_{[1:N]}$, $\mu^+_{[1:N-1]}$, $\mu^-_{[1:N-1]}$, $\epsilon_{\mathrm{abs}}$, $\epsilon_{\mathrm{rel}}$, $\epsilon_{\mathrm{buff}}$, $j_{\max}$
    \end{flushleft}
    \begin{algorithmic}[1]
    \Require $\hat{x}_k$, $\hat{u}_k$, $\mu^+_k$, $\mu^-_k$ are not all $0$
    \vspace{1em}
    \State $\sigma \leftarrow 0$\vspace{1ex}
    \For {$k=1,\ldots,N-1$}
    \State $\sigma \leftarrow \sigma + \norm{\hat{x}_{k}}_{2}^{2} + \norm{\hat{u}_{k}}_{2}^{2} + \norm{\mu^+_{k}}_{2}^{2} + \norm{\mu^-_{k}}_{2}^{2}$
    \EndFor\vspace{1ex}
    \State $\sigma \leftarrow \sigma + \norm{\hat{x}_{N}}_{2}^{2} + \norm{\hat{u}_{N}}_{2}^{2}$\vspace{1ex}
    \State $\sigma \leftarrow \sqrt{\sigma}$\vspace{1ex}
    \For {$j=1,\ldots,j_{\max}$}\vspace{1ex}
    \For {$k = 1,\ldots,N-1$}\vspace{1ex}
    \State $\phi_{k} \leftarrow \frac{1}{\sigma}\!\left(\hat{A}^{-}_{k}\,\hat{x}_{k} + \hat{A}^{+}_{k}\,\hat{x}_{k+1} + \hat{B}^{-}_{k} \hat{u}_{k} + \hat{B}^{+}_{k} \hat{u}_{k+1} + \mu^+_{k} - \mu^-_{k}\right)$\vspace{1ex}
    \State $\theta_{k} \leftarrow \frac{1}{\sigma} E_y\left(\hat{x}_{k+1} - \hat{x}_{k}\right)$\vspace{1ex}
    \EndFor\vspace{1ex}
    \State $\hat{x}_{1} \leftarrow \hat{A}^{-^{\top}}_{1}\!\phi_{1} - E_y^\top \theta_{1}$
    \State $\hat{u}_{1} \leftarrow \hat{B}^{-^{\top}}_{1}\!\phi_{1}$
    \State $\mu^+_{1} \leftarrow \phi_{1}$
    \State $\mu^-_{1} \leftarrow -\phi_{1}$\vspace{1ex}
    \For {$k = 2 ,\ldots,N-1$}
    \State $\hat{x}_{k} \leftarrow \hat{A}^{-^{\top}}_{k}\!\phi_{k} + \hat{A}^{+^{\top}}_{k-1} \phi_{k-1} - E_y^\top{v}_{k} + E_y^\top{v}_{k-1}$
    \State $\hat{u}_{k} \leftarrow \hat{B}^{-^{\top}}_{k}\!\phi_{k} + \hat{B}^{+^{\top}}_{k-1}\phi_{k-1}$
    \State $\mu^+_{k} \leftarrow \phi_{k}$
    \State $\mu^-_{k} \leftarrow -\phi_{k}$\vspace{0.5ex}
    \EndFor\vspace{1ex}
    \State $\hat{x}_{N} \leftarrow \hat{A}^{+^{\top}}_{N-1}\phi_{N-1} + E_y^\top{v}_{N-1}$
    \State $\hat{u}_{N} \leftarrow \hat{B}^{+^{\top}}_{N-1}\phi_{N-1}$\vspace{1ex}
    \State $\sigma^{\star} \leftarrow 0$\vspace{1ex}
    \For {$k =1,\ldots,N-1$}
    \State $\sigma^{\star} \leftarrow \sigma^{\star} + \norm{\hat{x}_{k}}_{2}^{2} + \norm{\hat{u}_{k}}_{2}^{2} + \norm{\mu^+_{k}}_{2}^{2} + \norm{\mu^-_{k}}_{2}^{2}$
    \EndFor\vspace{1ex}
    \State $\sigma^{\star} \leftarrow \sigma + \norm{\hat{x}_{N}}_{2}^{2} + \norm{\hat{u}_{N}}_{2}^{2}$\vspace{1ex}
    \State $\sigma^{\star} \leftarrow \sqrt{\sigma^{\star}}$\vspace{1ex}
    \If {$\abs{\sigma^{\star} - \sigma} \le \epsilon_{\mathrm{abs}} + \epsilon_{\mathrm{rel}}\,\max\{\sigma^{\star},\,\sigma\}$} \Comment{stopping criterion}
    \State \textbf{break}
    \ElsIf {$j < j_{\max}$}
    \State $\sigma \leftarrow \sigma^{\star}$
    \EndIf\vspace{1ex}
    \EndFor\vspace{1ex}
    \State $\sigma \leftarrow (1 + \epsilon_{\mathrm{buff}})\,\sigma^{\star}$\Comment{buffer the estimated maximum singular value}\vspace{0.125em}
    \end{algorithmic}
    \begin{flushleft}
        \textbf{Return:} $\sigma$
    \end{flushleft}
    \vspace{0.25em}
\end{algorithm}

\newgeometry{left=0.75in,
             right=0.75in,
             top=0.5in,
             bottom=0.5in,
             footskip=0in}
\begin{algorithm}[H]
\footnotesize
\caption{\pipgc }\label{alg:pipg_custom}
    \vspace{1ex}
    \begin{flushleft}
        \textbf{Inputs:} $w_{\mr{cost}}$, $e_{\mr{cost}}$, $w_{\mr{ep}}$, $w_{\mr{prox}}$\\ 
        \hphantom{\textbf{Inputs:}} $\hat{A}^{-}_{[1:N-1]}$, $\hat{A}^{+}_{[1:N-1]}$, $\hat{B}^{-}_{[1:N-1]}$, $\hat{B}^{+}_{[1:N-1]}$, $\hat{w}_{[1:N-1]}$\\
        \hphantom{\textbf{Inputs:}} $\hat{u}_{[1:N],\min}$, $\hat{u}_{[1:N],\max}$, $E_y$, $\hat{\varepsilon}_{[1:N-1]}$, $\tilde{E}_{\mr{i}}$, $\tilde{E}_{\mr{f}}$, $\hat{z}_{\mr{i}}$, $\hat{z}_{\mr{f}}$\\
        \hphantom{\textbf{Inputs:}} $\sigma$, $\omega$, $\rho$, $\epsilon_{\mathrm{abs}}$, $\epsilon_{\mathrm{rel}}$, $j_{\mathrm{check}}$, $j_{\max}$,\\
        \hphantom{\textbf{Inputs:}} $\hat{x}_{[1:N]}^{\star}$, $\hat{u}_{[1:N]}^{\star}$, ${\mu^{+,\star}}_{[1:N-1]}$, ${\mu^{-,\star}}_{[1:N-1]}$, $\phi^{\star}_{[1:N-1]}$, $\theta^{\star}_{[1:N-1]}$ \Comment{warm start}
    \end{flushleft}
    \vspace{-0.75em}
    \begin{algorithmic}[1]
    \State $ \tilde{x}_{[1:N]}^{1} \leftarrow  \hat{x}_{[1:N]}^{\star}$ \Comment{initialize primal variables}
    \State $ \tilde{u}_{[1:N]}^{1} \leftarrow  \hat{u}_{[1:N]}^{\star}$
    \State ${\nu^{+,1}}_{[1:N-1]} \leftarrow {\mu^{+,\star}}_{[1:N-1]}$
    \State ${\nu^{-,1}}_{[1:N-1]} \leftarrow {\mu^{-,\star}}_{[1:N-1]}$
    \State $\varphi_{[1:N-1]}^{1} \leftarrow \phi_{[1:N-1]}^{\star}$\Comment{initialize dual variables}
    \State $\vartheta_{[1:N-1]}^{1} \leftarrow \theta_{[1:N-1]}^{\star}$\vspace{1ex}
    \State $\alpha \leftarrow \frac{2}{w_{\mr{prox}} + \sqrt{w_{\mr{prox}}^{2} + 4\omega\sigma}}$\Comment{step-sizes} 
    \State $\beta \leftarrow \omega\alpha$\vspace{1ex}
    \For {$j = 1,\ldots,j_{\max}$}\vspace{1ex}
    \State $ \hat{x}^{j}_{1} \leftarrow  \tilde{x}_{1}^{j}-\alpha\,(w_{\mr{prox}}\, \tilde{x}_{{1}}^{j} + \hat{A}^{-\top}_{1}\!\varphi_{1}^{j} - E_y^\top\vartheta_{1}^{j})$
    \State $\tilde{E}_{\mr{i}}\hat{x}^{j}_1\gets \hat{z}_{\mr{i}}$
    \State $ \hat{u}^{j}_{1} \leftarrow \max\{\hat{u}_{1,\min},\min\{\hat{u}_{1,\max},\tilde{u}_{{1}}^{j}-\alpha\,(w_{\mr{prox}}\, \tilde{u}_{{1}}^{j} + \hat{B}^{-\top}_{1}\!\varphi^{j}_{1})\}\}$\vspace{1ex}
    \For {$k = 2,\ldots,N-1$}\Comment{projected gradient step}\vspace{1ex}
    \State $ \hat{x}_{k}^{j+1} \leftarrow \tilde{x}_{{k}}^{j}-\alpha\,(w_{\mr{prox}}\, \tilde{x}_{{k}}^{j}  + \hat{A}^{-\top}_{k}\!\varphi_{k}^{j} + \hat{A}^{+\top}_{k-1}\varphi_{k-1}^{j} - E_y^\top\vartheta_{k}^{j} + E_y^\top\vartheta_{k-1}^{j})$
    \State $ \hat{u}_{k}^{j+1} \leftarrow  \max\{\hat{u}_{k,\min},\min\{\hat{u}_{k,\max},\tilde{u}_{{k}}^{j}-\alpha\,(w_{\mr{prox}}\, \tilde{u}_{{k}}^{j} + \hat{B}^{-\top}_{k}\!\varphi_{k}^{j} + \hat{B}^{+{\top}}_{k-1}\vartheta_{k-1}^{j})\}\}$\vspace{1ex}
    \EndFor\vspace{1ex}
    \State $ \hat{x}_{N}^{j+1} \leftarrow \tilde{x}_{{N}}^{j} - \alpha\,(w_{\mr{prox}}\, \tilde{x}_{{N}}^{j} + w_{\mr{cost}}e_{\mr{cost}} + \hat{A}^{+{\top}}_{N-1}\varphi_{N-1}^{j} + E_y^\top\vartheta_{N-1}^{j})$
    \State $\tilde{E}_{\mr{f}}\hat{x}_N^{j+1}\gets \hat{z}_{\mr{f}}$
    \State $ \hat{u}_{N}^{j+1} \leftarrow \max\{\hat{u}_{N,\min},\min\{\hat{u}_{N,\max},\tilde{u}_{{N}}^{j} - \alpha\,(w_{\mr{prox}}\, \tilde{u}_{{N}}^{j}  + \hat{B}^{+{\top}}_{N-1}\varphi_{N-1}^{j})\}\}$
    \State $\mu^{\pm,j}_{[1:N-1]} \leftarrow \max\{0,\nu^{\pm,j}_{{[1:N-1]}}-\alpha\,(w_{\mr{ep}}1_{n_x} + \varphi^{j}_{[1:N-1]})\}$\vspace{1ex}
    \For {$k = 1,\ldots,N-1$}\Comment{PI feedback of affine equality constraint violation}\vspace{1ex}
    \State $\phi_{k}^{j+1} \leftarrow \varphi_{k}^{j} + \beta\,(\hat{A}^{-}_{k}\,(2 \hat{x}^{j+1}_{k} -  \tilde{x}^{j}_{{k}}) + \hat{A}^{+}_{k}\,(2 \hat{x}^{j+1}_{k+1} -  \tilde{x}^{j}_{{k+1}}) + \hat{B}^{-}_{k}\,(2 \hat{u}^{j+1}_{k} -  \tilde{u}^{j}_{{k}}) + \hat{B}^{+}_{k}\,(2 \hat{u}^{j+1}_{k+1} -  \tilde{u}^{j}_{{k+1}})$
    \State \hphantom{$\hat{w}_{k}^{j} \leftarrow \eta_{k}^{j} + \alpha\,($}$+\:(2\,\mu^{+,j+1}_{k} - \nu^{+,j}_{{k}}) - (2\,\mu^{-,j+1}_{k} - \nu^{-,j}_{{k}}) + \hat{w}_{k})$\vspace{1ex}
    \State $\theta_{k}^{j+1} \leftarrow \max\{0,\, \vartheta_{k}^{j} + \beta\,((2 E_y\hat{x}^{j+1}_{k+1} - E_y\tilde{x}^{j}_{{k+1}}) - (2E_y\hat{x}^{j+1}_{k} -  E_y\tilde{x}^{j}_{{k}})  - \hat{\varepsilon}_k)\}$\vspace{1ex}
    \EndFor\vspace{1ex}
    \State $ \tilde{x}_{{[1:N]}}^{j+1} \leftarrow (1 - \rho)\, \tilde{x}_{{[1:N]}}^{j} + \rho\, \hat{x}_{[1:N]}^{j+1}$ \Comment{extrapolate primal variables}
    \State $ \tilde{u}_{{[1:N]}}^{j+1} \leftarrow (1 - \rho)\, \tilde{u}_{{[1:N]}}^{j} + \rho\, \hat{u}_{[1:N]}^{j+1}$
    \State $ \nu^{\pm,j+1}_{{[1:N-1]}} \leftarrow (1 - \rho)\,\nu^{\pm,j}_{{[1:N-1]}} + \rho\,\mu^{\pm,j+1}_{[1:N-1]}$
    \State $ \varphi_{[1:N-1]}^{j+1} \leftarrow (1 - \rho)\,\varphi_{[1:N-1]}^{j} + \rho\,\phi_{[1:N-1]}^{j+1}$\Comment{extrapolate dual variables}
    \State $\vartheta_{[1:N-1]}^{j+1} \leftarrow (1 - \rho)\,\vartheta_{[1:N-1]}^{j} + \rho\,\theta_{[1:N-1]}^{j+1}$\vspace{1ex}
    \If {$j \operatorname{mod} j_{\mathrm{check}} = 0$} \Comment{check stopping criterion every $j_{\mathrm{check}}$ iterations}\vspace{1ex}
    \State $\textsc{terminate} \leftarrow \textsc{stopping}_{\text{custom}}(\hat{x}_{[1:N]}^{j+1},\,\hat{u}_{[1:N]}^{j+1},\,\mu_{[1:N-1]}^{\pm,j+1},\,\phi_{[1:N-1]}^{j+1},\,\theta_{[1:N-1]}^{j+1}$
    \State \hphantom{$\textsc{terminate} \leftarrow \textsc{stopping}_{\text{custom}}\:$}$\,\hat{x}_{[1:N]}^{j},\,\hat{u}_{[1:N]}^{j},\,\mu_{[1:N-1]}^{\pm,j},\,\phi_{[1:N-1]}^{j},\,\theta_{[1:N-1]}^{j},\,\epsilon_{\mathrm{abs}},\,\epsilon_{\mathrm{rel}})$\vspace{1ex}
    \If {$\textsc{terminate} = \textsc{true}$}\vspace{1ex}\Comment{stopping criterion}
    \State \textbf{break}\vspace{1ex}
    \EndIf\vspace{1ex}
    \EndIf\vspace{1ex}
    \EndFor\vspace{1ex}
    \State $ \hat{x}_{[1:N]}^{\star} \leftarrow  \hat{x}_{[1:N]}^{j+1}$ \Comment{update primal variables}
    \State $ \hat{u}_{[1:N]}^{\star} \leftarrow  \hat{u}_{[1:N]}^{j+1}$
    \State $\mu_{[1:N-1]}^{\pm,\star} \leftarrow \mu_{[1:N-1]}^{\pm,j+1}$
    \State $\phi_{[1:N-1]}^{\star} \leftarrow \phi_{[1:N-1]}^{j+1}$ \Comment{update dual variables}
    \State $\theta_{[1:N-1]}^{\star} \leftarrow \theta_{[1:N-1]}^{j+1}$
    \end{algorithmic}
    \begin{flushleft}
        \textbf{Return:} $ \hat{x}_{[1:N]}^{\star}$, $ \hat{u}_{[1:N]}^{\star}$, $\mu^{\pm,\star}_{[1:N-1]}$,  $\phi^{\star}_{[1:N-1]}$, $\theta^{\star}_{[1:N-1]}$
    \end{flushleft}
\end{algorithm}

\restoregeometry

\begin{algorithm}[H]
\small
\caption{Customized stopping criterion evaluation:\\[1ex]
\hphantom{\textbf{Algorithm 7}~}$\textsc{stopping}_{\text{custom}}(\hat{x}_{[1:N]}^{j+1},\,\hat{u}_{[1:N]}^{j+1},\,\mu_{[1:N-1]}^{+,j+1},\,\mu_{[1:N-1]}^{-,j+1},\,\phi_{[1:N-1]}^{j+1},\,\theta_{[1:N-1]}^{j+1}$\\
\hphantom{\textbf{Algorithm 7}~$\textsc{stopping}_{\text{custom}}$\:}$\,\hat{x}_{[1:N]}^{j},\,\hat{u}_{[1:N]}^{j},\,\mu_{[1:N-1]}^{+,j},\,\mu_{[1:N-1]}^{-,j},\,\phi_{[1:N-1]}^{j},\,\theta_{[1:N-1]}^{j},\,\epsilon_{\mathrm{abs}},\,\epsilon_{\mathrm{rel}})$}\label{alg:stopping_custom}
    \vspace{0.25em}
    \begin{flushleft}
        \textbf{Inputs:} $\hat{x}_{[1:N]}^{j+1},\,\hat{u}_{[1:N]}^{j+1},\,\mu_{[1:N-1]}^{+,j+1},\,\mu_{[1:N-1]}^{-,j+1},\,\phi_{[1:N-1]}^{j+1},\,\theta_{[1:N-1]}^{j+1}$\\
        \hphantom{\textbf{Inputs:}} $\,\hat{x}_{[1:N]}^{j},\,\hat{u}_{[1:N]}^{j},\,\mu_{[1:N-1]}^{+,j},\,\mu_{[1:N-1]}^{-,j},\,\phi_{[1:N-1]}^{j},\,\theta_{[1:N-1]}^{j},\,\epsilon_{\mathrm{abs}},\,\epsilon_{\mathrm{rel}}$
    \end{flushleft}
    \begin{algorithmic}[1]
    \State $z^{j+1}_{\infty} \leftarrow \max\!\left\{\|\hat{x}_{[1:N]}^{j+1}\|_{\infty},\, \|\hat{u}_{[1:N]}^{j+1}\|_{\infty},\, \|\mu_{[1:N-1]}^{+,j+1}\|_{\infty},\, \|\mu_{[1:N-1]}^{-,j+1}\|_{\infty}\right\}$
    \State $z^{j\hphantom{+ 1}}_{\infty} \leftarrow \max\!\left\{\|\hat{x}_{[1:N]}^{j}\|_{\infty},\, \|\hat{u}_{[1:N]}^{j}\|_{\infty},\, \|\mu_{[1:N-1]}^{+,j}\|_{\infty},\, \|\mu_{[1:N-1]}^{-,j}\|_{\infty}\right\}$
    \State $z^{\Delta j}_{\infty} \,\,\leftarrow \max\!\left\{\|\hat{x}_{[1:N]}^{j+1} - \hat{x}_{[1:N]}^{j}\|_{\infty},\, \|\hat{u}_{[1:N]}^{j+1} - \hat{u}_{[1:N]}^{j}\|_{\infty},\, \|\mu_{[1:N-1]}^{+,j+1} - \mu_{[1:N-1]}^{+,j}\|_{\infty},\, \|\mu_{[1:N-1]}^{-,j+1} - \mu_{[1:N-1]}^{-,j}\|_{\infty}\right\}$\vspace{1ex}
    \State $r^{j+1}_{\infty} \leftarrow \max\!\left\{\|\phi_{[1:N-1]}^{j+1}\|_{\infty},\, \|\theta_{[1:N-1]}^{j+1}\|_{\infty}\right\}$
    \State $r^{j\hphantom{+ 1}}_{\infty} \leftarrow \max\!\left\{\|\phi_{[1:N-1]}^{j}\|_{\infty},\, \|\theta_{[1:N-1]}^{j}\|_{\infty}\right\}$
    \State $r^{\Delta j}_{\infty} \,\, \leftarrow \max\!\left\{\|\phi_{[1:N-1]}^{j+1} - \phi_{[1:N-1]}^{j}\|_{\infty},\, \|\theta_{[1:N-1]}^{j+1} - \theta_{[1:N-1]}^{j}\|_{\infty}\right\}$\vspace{1ex}
    \If {$z^{\Delta j}_{\infty} \leq \epsilon_{\text{abs}} + \epsilon_{\text{rel}}\,\max\!\left\{z^{j+1}_{\infty},\,z^{j}_{\infty}\right\}$ \textbf{and} $r^{\Delta j}_{\infty} \leq \epsilon_{\text{abs}} + \epsilon_{\text{rel}}\,\max\!\left\{r^{j+1}_{\infty},\,r^{j}_{\infty}\right\}$}\vspace{1ex}
    \State $\textsc{terminate} \leftarrow \textsc{true}$\vspace{1ex}
    \Else
    \State $\textsc{terminate} \leftarrow \textsc{false}$\vspace{1ex}
    \EndIf\vspace{1ex}
    \end{algorithmic}
    \begin{flushleft}
        \textbf{Return:} \textsc{terminate}
    \end{flushleft}
    \vspace{0.25em}
\end{algorithm}

%% file: sections/5-pdg-formulation.tex
\section{6-DoF Powered Descent Guidance}\label{sec:problem-formulation}

The 6-DoF rocket model is based on \cite{szmuk2020successive}, with the inclusion of an independent torque control input by means of reaction control thrusters (as shown in \cite{kamath2023customized}), and with the aerodynamics forces and moments ignored.

The system dynamics \eqref{sys-dyn} and constraint function \eqref{path-cnstr} are given by:
\begin{subequations}
\begin{align}
    & \xi = (m,\bm{r}_{\mc{I}},\bm{v}_{\mc{I}},\bm{q}_{\mc{I}\leftarrow\mc{B}},\bm{\omega}_{\mc{B}})\\
    & \zeta = (\bm{T}_{\mc{B}}, \bm{\gamma}_{\mc{B}})\\
    & F(\xi,\zeta) = \left[ \begin{array}{c} -\alpha_{\dot{m}}\|\bm{T}_{\mc{B}}\|_2\\
                                             \bm{v}_{\mc{I}}\\
                                             \frac{1}{m}C_{\mc{I}\leftarrow\mc{B}}\bm{T}_{\mc{B}}+\bm{g}_{\mc{I}}\\
                                             \frac{1}{2}\Omega(\bm{\omega}_{\mc{B}})\bm{q}_{\mc{I}\leftarrow\mc{B}}\\
                                             J_{\mc{B}}^{-1}\left( \bm{r}_{T,\mc{B}}\times \bm{T}_{\mc{B}} - \bm{\omega}_{\mc{B}}\times J_{\mc{B}}\bm{\omega}_{\mc{B}} + \bm{\gamma}_{\mc{B}} \right)
    \end{array}\right]\\
    & g(\xi,\zeta) = \left[\begin{array}{c}
         m_{\mr{dry}} - m  \\
         -\bm{e}_1^\top\bm{r}_{\mc{I}}\\
         \|\bm{v}_{\mc{I}}\|_2^2 - v_{\max}^2\\
         4\|H_\theta\bm{q}_{\mc{I}\leftarrow\mc{B}}\|_2^2 + (1-\cos\theta_{\max})^2\\
         \|\bm{\omega}_{\mc{B}}\|_2^2 - \omega_{\max}^2\\
         \|\bm{T}_{\mc{B}}\|_2 - \bm{e}_1^\top \bm{T}_{\mc{B}}\sec\delta_{\max}\\
         \|\bm{T}_{\mc{B}}\|_2 - T_{\max}\\
        -\|\bm{T}_{\mc{B}}\|_2 + T_{\min}\\
         \|\bm{\gamma}_{\mc{B}}\|^2_2 - \gamma_{\max}^2 
    \end{array}\right]
\end{align}
\end{subequations}

The numerical values for the parameters within the system dynamics $(\alpha_{\dot{m}},\bm{g}_{\mc{I}},J_{\mc{B}},\bm{r}_{T,\mc{B}},H_{\theta})$, boundary conditions ($z_{\mr{i}},z_{\mr{f}}$), and constraint bounds $(m_{\mr{dry}},\gamma_{\mr{gs}},v_{\max},\theta_{\max},T_{\max},T_{\min})$ correspond to the baseline problem in \cite{miki2020successive}. Note that the non-differentiability of the constraint functions for thrust magnitude bounds at $\bm{T}_{\mc{B}} = 0_3$ does not affect the linearization operations within SCP in practice with an appropriate initial guess for the thrust vector.
Furthermore, the box constraint on the scaled control input in Lines 12, 15 and 19 of Algorithm \ref{alg:pipg_custom} are only imposed on the dilation factor.

%% file: sections/6-gpu.tex
\section{GPU Architecture and Code Implementation}\label{section:gpu-architecture}

\subsection{Overview}

Loosely speaking, CPUs are designed to execute code serially. Consider the problem of adding two vectors. A CPU is designed to add these two vectors by looping through the elements of each vector and adding the corresponding elements. In contrast, GPUs have many specialized hardware units that allow them to perform many floating point operations in parallel. A GPU would add two vectors together by adding the corresponding elements in parallel by assigning each scalar addition to a different compute thread.

In GPU programming, a single sequence of instructions is called a \textit{thread}. Back to the example of adding two vectors on a GPU, a single thread would execute a single scalar addition. Threads in turn are organized into a \textit{grid} of \textit{thread-blocks}. In the vector addition example, we would want to have the same number of threads as the number of elements in each vector. To write programs for the GPU, NVIDIA provides the compute unified device architecture ({\cuda}) platform, which is an extension of the C programming language \cite{cudaGuide}. 

{\cuda} enables the programmer to write special GPU functions called \textit{kernels}, which, when called, are executed a specified number of times in parallel. When the programmer writes a kernel, they are specifying how a single thread should execute the given code; when calling a kernel, the programmer specifies the number of threads per thread-block and the number of thread-blocks in the grid. The product of these two numbers is the number of threads being called. Each thread of execution has a unique identifier based on its thread-block index within the grid and its thread index within the thread-block. The thread ID in a kernel is an integer in $[0, N)$, where $N$ is the number of threads called. However, before performing computations on the GPU, we must first ensure that the data the kernel is operating on is in GPU memory rather than in CPU memory. This typically requires the programmer to allocate memory on the GPU to store this data, using the \texttt{cudaMalloc()} function, and to then copy the data from the CPU to the GPU using the \texttt{cudaMemcpy()} function.

When the kernel is called, the thread-blocks are then distributed to \textit{streaming multiprocessors} (SM), which contain a collection of {\cuda} cores. These {\cuda} cores execute the threads within a thread-block. When a SM finishes executing all threads on a thread-block, the GPU allocates it another thread-block.

In all of this, the programmer only has to worry about creating a grid of thread-blocks when they call the kernel. The GPU is responsible for allocating thread-blocks to SMs, which are in turn are responsible for allocating threads to {\cuda} cores. The thread hierarchy is summarized by Figure \ref{fig:gpu-compute}.

As an example, consider the code in Listing \ref{code:cudasnippet}. The first function is how one would write a function to add two vectors on a CPU. The second function is how one would write a {\cuda} kernel to add two vectors in parallel. When launching this kernel, each GPU thread executes the code within the kernel in parallel. The first line in this kernel computes the unique thread ID for the thread; subsequently this thread will perform one scalar addition for the index that matches its thread ID. The \texttt{main()} function shows how to call the serial version of this vector addition and how to call the {\cuda} kernel. We must first allocate the memory for arrays $a$, $b$, and $c$ on the GPU and then copy the data in the arrays from the CPU to the GPU. Finally, the kernel can be called with the desired grid size and thread-block size.

\input{tikz/gpu-model}

\begin{minipage}{\linewidth}

\begin{lstlisting}[language=C++, caption=Code for serial (CPU) and parallelized (GPU) vector addition, label={code:cudasnippet}]
/* Performs the vector addition c = a + b serially, where all vectors have N elements. */
void add_vector_serial(int* a, int* b, int* c, int N)
{
    for(int i = 0; i < N; i++)
    {
        c[i] = a[i] + b[i];
    }
}

/* CUDA kernel to perform c = a + b, where all vectors have N elements. */
__global__ void add_vector_parallel(int* a, int* b, int* c, int N)
{
    /* Get the ID of the thread executing this function. */
    int thread_id = blockIdx.x*blockDim.x + threadIdx.x;

    /* Safeguard if the kernel is launched with more threads than the number of elements. */
    if (i < N)
    {
        c[i] = a[i] + b[i];
    }
}

int main()
{
    /* Declare pointers for arrays in CPU memory. */
    int* a[N];
    int* b[N];
    int* c[N];

    /* Helper function to fill a and b with data. */
    initialize(a, b);

    /* Declare pointers for arrays in GPU memory. */
    int* d_a, d_b, d_c;

    /* Allocate memory for arrays in GPU memory. */
    cudaMalloc((void**) &d_a, N * sizeof(int));
    cudaMalloc((void**) &d_b, N * sizeof(int)));
    cudaMalloc((void**) &d_c, N * sizeof(int)));

    /* Copy data in array from CPU to GPU memory. */
    cudaMemcpy(d_a, a, N * sizeof(int), cudaMemcpyHostToDevice);
    cudaMemcpy(d_b, b, N * sizeof(int), cudaMemcpyHostToDevice);
    cudaMemcpy(d_c, c, N * sizeof(int), cudaMemcpyHostToDevice);
    
    /* Call serial addition function on CPU. */
    add_vector_serial(a, b, c, N);

    /* Call parallel addition kernel on GPU. */
    add_vector_parallel<<<grid_size, block_size>>>(d_a, d_b, d_c, N);
}

\end{lstlisting}
\end{minipage}

\subsection{Code Structure}
The data transfer bandwidth between the CPU and the GPU is much slower than that between the {\cuda} cores and the GPU memory. Thus to achieve peak performance, it is important to minimize data transfer between the CPU and GPU \cite{cudaGuide}. Consequently, we chose the following execution model: create an array of problem instances with varying initial conditions on the CPU, allocate the necessary memory on the GPU to store these problem instances, copy the problem instance data to the GPU, launch a kernel that solves all of these trajectory optimization problem instances on the GPU, and finally copy back the trajectories from the GPU to the CPU. In this paradigm, we are assigning each thread one trajectory to compute. Each thread is able to access the problem instance it needs to solve by indexing into the array of problem data using its unique thread ID. Figure \ref{fig:data-flow} shows the flow of data in our implementation.

By taking this approach we only have to incur the costly overhead of CPU-GPU communication once when sending the problem instances to the GPU and once when copying the computed trajectories back from the GPU.

As mentioned in Section \ref{sec:intro}, a key difficulty with implementing this procedure on the GPU is the need to implement all operations from scratch in {\cuda}, since we are unable to use off-the-shelf linear algebra routines and convex optimization solvers. Therefore, all code necessary to solve for a trajectory was implemented in {\cuda}, including Level 1, 2, and 3 BLAS operations and the {\pipgc} algorithm.

\input{tikz/data-flow}

%% file: tikz/gpu-model.tex
\begin{figure}[H]
\begin{mybox}
\centering

\begin{tikzpicture}
\draw[->,decorate,decoration={snake,amplitude=2pt,segment length=5pt}] (0,0) -- (0,-0.7);
\node[above, font=\scriptsize] at (0,0.1) {Thread};
\draw[->] (1,-0.35) -- (4,-0.35) node[midway, above, font=\scriptsize] {executed by};

\node[draw, rectangle, minimum size=1cm] at (5,-0.35) {\large $+$};
\node[above, font=\scriptsize] at (5,0.1) {{\cuda} Core};

\foreach \i in {1,...,4} {
\draw[->,decorate,decoration={snake,amplitude=2pt,segment length=5pt}] (0.15*\i-0.5,-2.3) -- (0.15*\i-0.5,-3);
}
\node[draw=none] (ellipsis1) at (0.45,-2.65) {$\cdots$};
\draw[->,decorate,decoration={snake,amplitude=2pt,segment length=5pt}] (1.25-0.5,-2.3) -- (1.25-0.5,-3);
\draw[draw=black, thick] (-0.5,-3.1) rectangle (0.9,-2.2);
\node[above, font=\scriptsize] at (0.25,-2.2) {Thread-Block};
\draw[->] (1,-2.55) -- (4,-2.55) node[midway, above, font=\scriptsize] {executed by};

\draw[draw=black, thick, fill=white] (4.5,-3.25) rectangle (8.75,-1.75);
\foreach \i in {1,...,2} {
\node[draw, minimum size=1cm] at (1.1*\i+4.15,-2.5) {\large $+$};
}
\node[above, font=\scriptsize] at (6.5,-1.75) {Streaming Multiprocessor};
\node[draw=none] (ellipsis1) at (7.2,-2.5) {$\cdots$};
\node[draw, minimum size=1cm] at (8,-2.5) {\large $+$};

\foreach \j in {-2, ..., 0}
{
\foreach \i in {1,...,4} {
    \draw[->,decorate,decoration={snake,amplitude=2pt,segment length=5pt}] (0.15*\i-0.8+1.5*\j,-4.3) -- (0.15*\i-0.8+1.5*\j,-5);
  }
  \node[draw=none] (ellipsis1) at (0.15+1.5*\j,-4.65) {$\cdots$};
  \draw[->,decorate,decoration={snake,amplitude=2pt,segment length=5pt}] (1.25-0.8+1.5*\j,-4.3) -- (1.25-0.8+1.5*\j,-5);
  \draw[draw=black, thick] (-0.8+1.5*\j,-5.1) rectangle (0.6+1.5*\j,-4.2);
}
  \node[above, font=\scriptsize] at (-1.5,-4) {Grid};
  \draw[->] (1,-4.55) -- (4,-4.55) node[midway, above, font=\scriptsize] {executed by};
\draw[draw=black, thick] (-4.05,-5.3) rectangle (0.9,-4);

\foreach \j in {-2,..., 0}
{
  \draw[draw=black, thick, fill=white] (4.5+0.1*\j,-5.5-0.1*\j) rectangle (8.75+0.1*\j,-4-0.1*\j);
  \foreach \i in {1,...,2} {
    \node[draw, minimum size=1cm] at (1.1*\i+4.15,-4.75) {\large $+$};
  }
  \node[above, font=\scriptsize] at (6.5,-3.85) {GPU};
  \node[draw=none] (ellipsis1) at (7.2,-4.75) {$\cdots$};
  \node[draw, minimum size=1cm] at (8,-4.75) {\large $+$};
}
\end{tikzpicture}
\end{mybox}
\caption{Thread hierarchy in the GPU}
\label{fig:gpu-compute}
\end{figure}
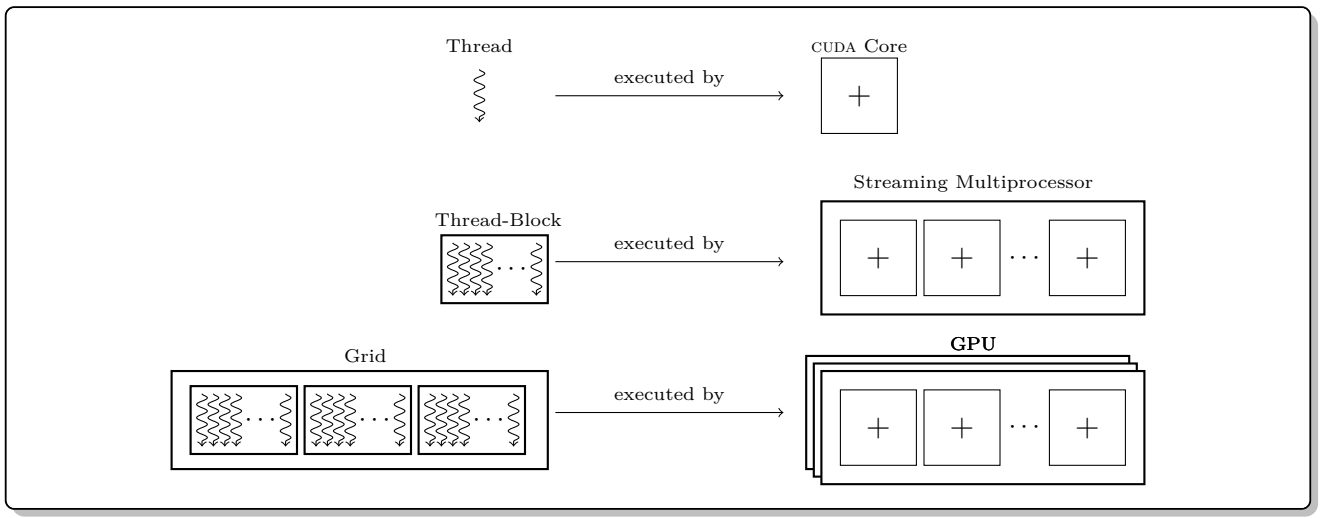

%% file: tikz/data-flow.tex
\begin{figure}
\begin{mybox}
\centering
\begin{tikzpicture}[scale=0.8, transform shape]

\draw (0,2) rectangle node {Initial Conditions} ++(3,1);

\draw[->] (1.5,2) -- node[above] {} ++(0,-1);

\draw (0,0) rectangle node {CPU Memory} ++(3,1);

\draw[->] (3,0.5) -- node[above] {} ++(1,0);

\draw (4,0) rectangle node {GPU Memory} ++(2.5,1);

\draw[->] (6.5,0.5) -- node[above] {} ++(1,0);


\draw (7.5,-1.5) -- (7.5,3);

\foreach \y in {0, 2, 3}
  \draw[->] (7.5,1.5*\y-1.5) -- ++(1,0);

\foreach \y in {0, 2, 3}
    \draw (8.5,-2+1.5*\y) rectangle node {GPU Thread} ++(2.5,1);

\node[draw=none] (ellipsis) at (9.75,0) [font=\Huge] {$\vdots$};

\foreach \y in {0, 2, 3}
  \draw[->] (11,1.5*\y-1.5) -- ++(1,0);

  \draw (12,-1.5) -- (12,3);

  \draw[->] (12,0.5) -- node[above] {} ++(1,0);

\draw (13,0) rectangle node {CPU Memory} ++(2.5,1);

  \draw[->] (15.5,0.5) -- node[above] {} ++(1,0);

\draw (16.5,0) rectangle node {Output File} ++(2.5,1);

\end{tikzpicture}
\end{mybox}
\caption{Data flow between the CPU and the GPU}
\label{fig:data-flow}
\end{figure}
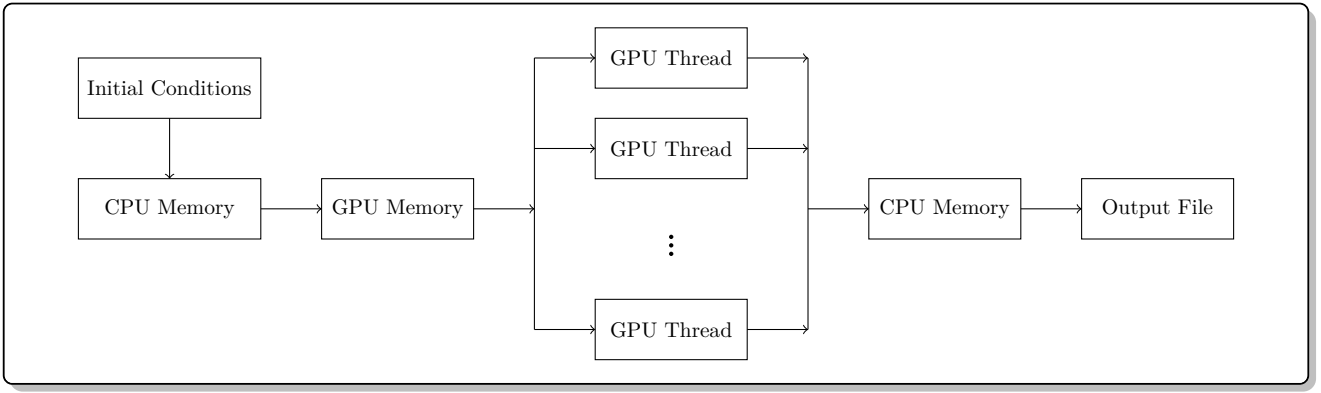

%% file: sections/7-numerical-results.tex
\section{Numerical Results}
This section presents and analyzes Monte Carlo results for the 6-DoF powered-descent guidance problem outlined in Section \ref{sec:problem-formulation} using the solution procedure in Section \ref{sec:ctcs}. We then compare the computation time between a serial CPU implementation and parallel GPU implementation for Monte Carlo simulations of various sizes. All results were obtained on a machine with an AMD Ryzen 9 7950X3D CPU and an NVIDIA RTX 3090 GPU.

When conducting the Monte Carlo, we disperse the initial position according to the following uniform distribution:

\begin{equation}
    \bm{r}_{\mc{I}}(0) \sim \begin{bmatrix}
        \mathcal{U}(6, 9)\\
        \mathcal{U}(3, 6)\\
        \mathcal{U}(1, 2) \\
    \end{bmatrix}
\end{equation}

\begin{figure}[ht]

        \begin{subfigure}{.5\textwidth}
            \centering
            \includegraphics[width=1\textwidth]{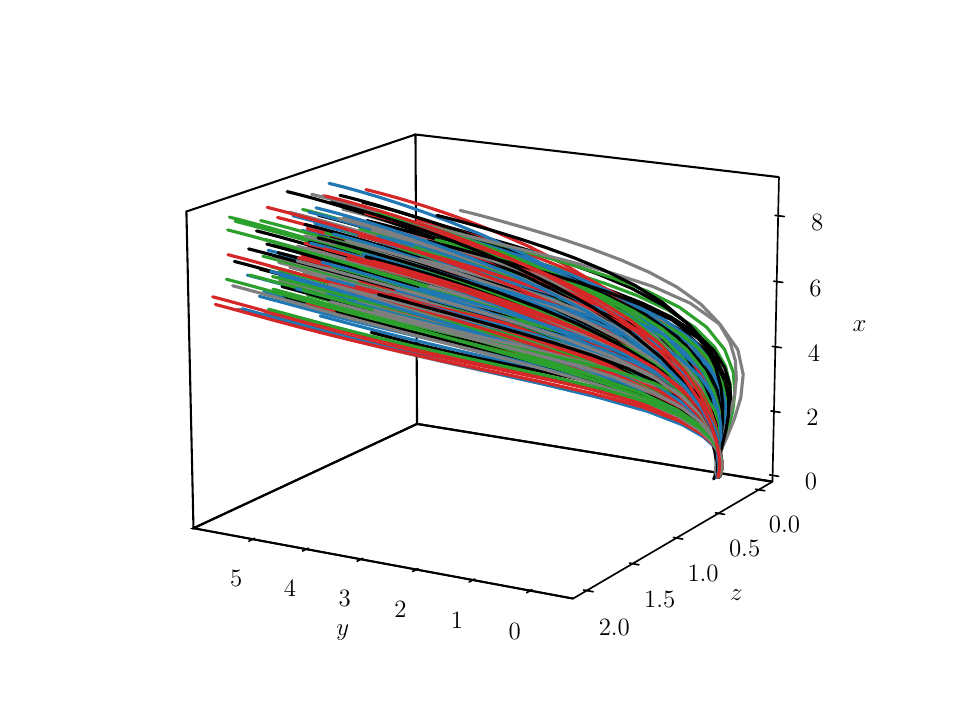}
        \end{subfigure}
        \begin{subfigure}{.5\textwidth}
            \centering
            \includegraphics[width=\textwidth]{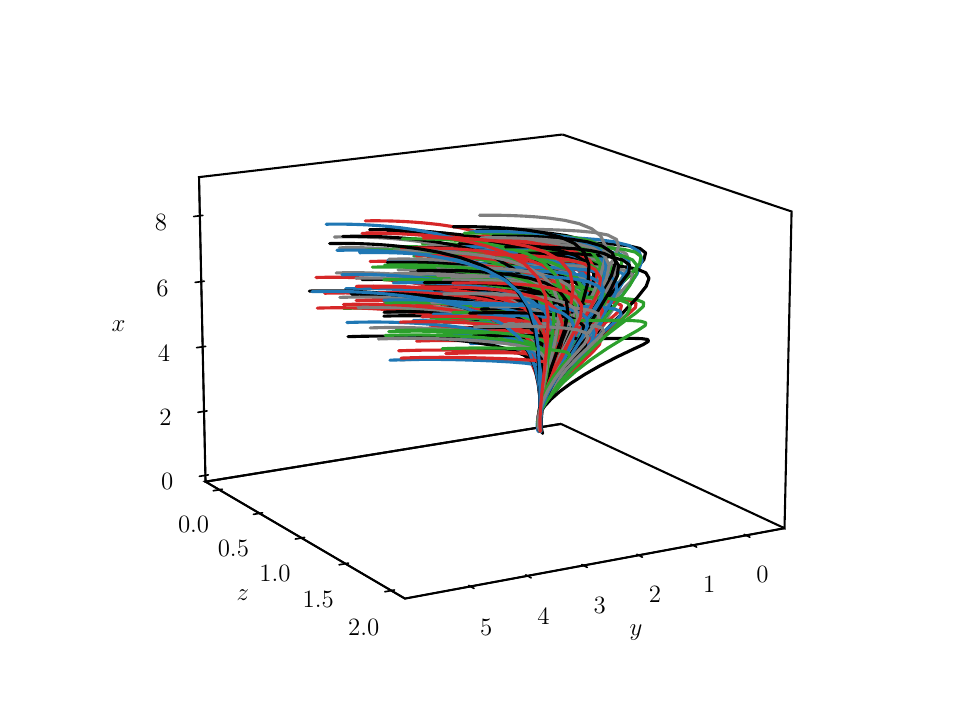}
        \end{subfigure}
        \caption{Trajectories generated from 256 Monte Carlo runs}
        \label{fig:mc-trajs}
        \end{figure}

\begin{figure}[ht]
        \begin{subfigure}{.5\textwidth}
            \centering
            \includegraphics[width=1\linewidth]{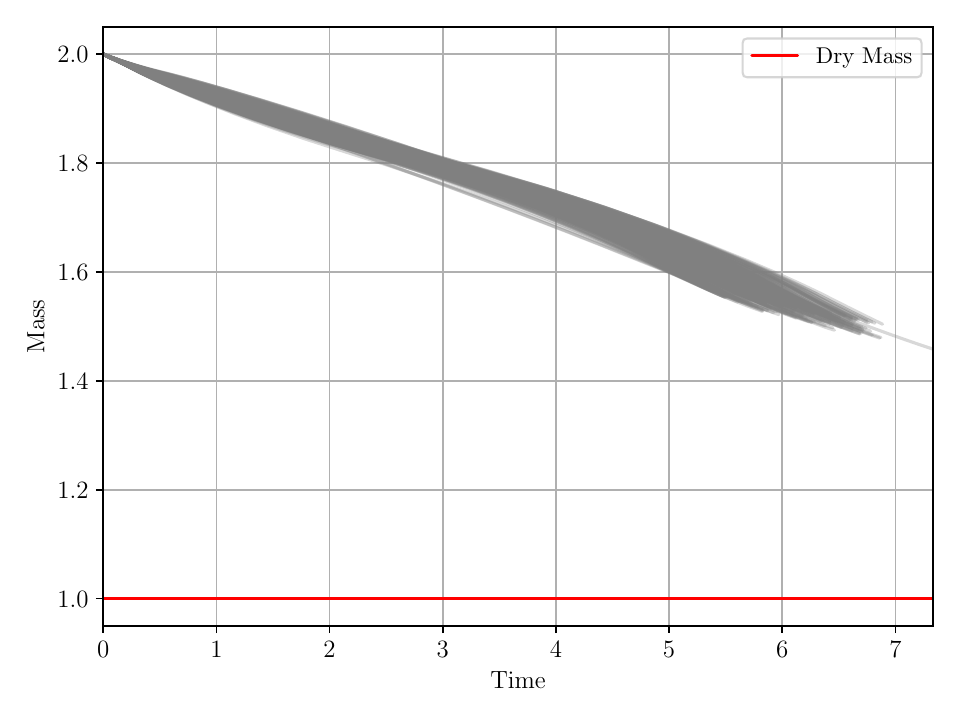}
            \caption{Mass Depletion}
        \end{subfigure}
        \begin{subfigure}{.5\textwidth}
            \centering
            \includegraphics[width=\linewidth]{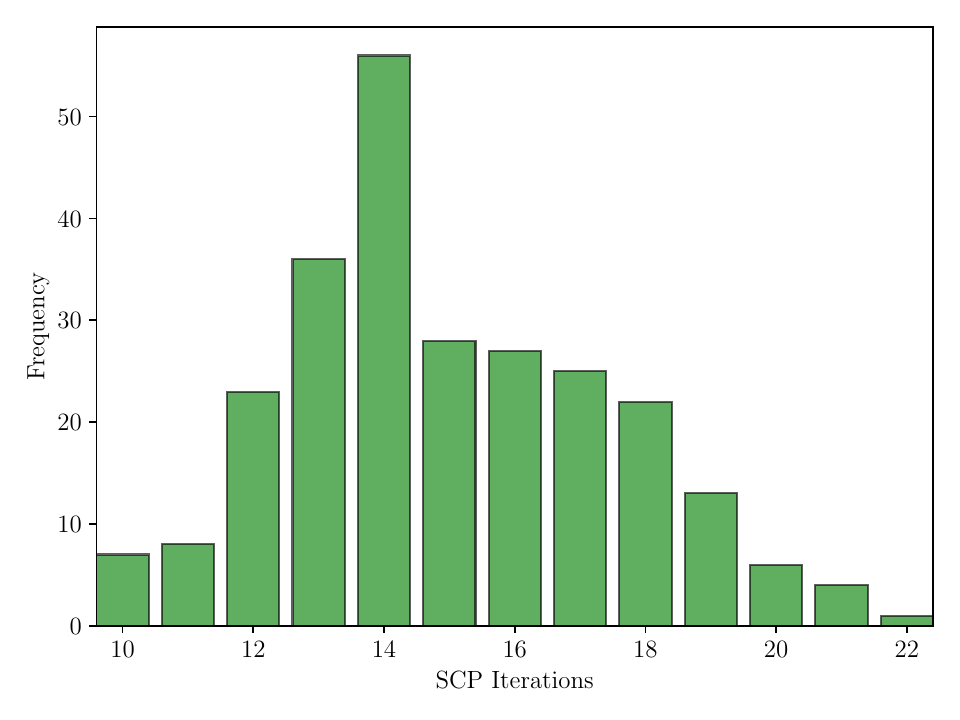}
            \caption{Distribution of SCP Iterations}
        \end{subfigure}
        \caption{Statistics from 256 Monte Carlo runs}
        \label{fig:mc-plots}
    \end{figure}

For all simulations, we use the same weights, $w_{\mr{cost}}$, $w_{\mr{tr}}$, and $w_{\mr{vc}}$, defined in Equation \ref{cvx-subproblem-custom:obj}, we set the maximum number of SCP iterations to $25$, and we use $2500$ {\pipg} iterations with the same hyperparameters $\omega$ and $\rho$ for each subproblem solve. We find that SCP converges within $25$ iterations for all $256$ of these trajectories, which is the first demonstration of robustness to variations in initial conditions for SCP with {\pipg} as a subproblem solver. Figures \ref{fig:mc-trajs} and \ref{fig:mc-plots} show the trajectories obtained. An extensive Monte Carlo analysis of the 6-DoF PDG problem is provided in \cite{malyuta2019discretization}.

    \begin{figure}[ht]

        \begin{subfigure}{.5\textwidth} 
            \centering
            \includegraphics[width=1\textwidth]{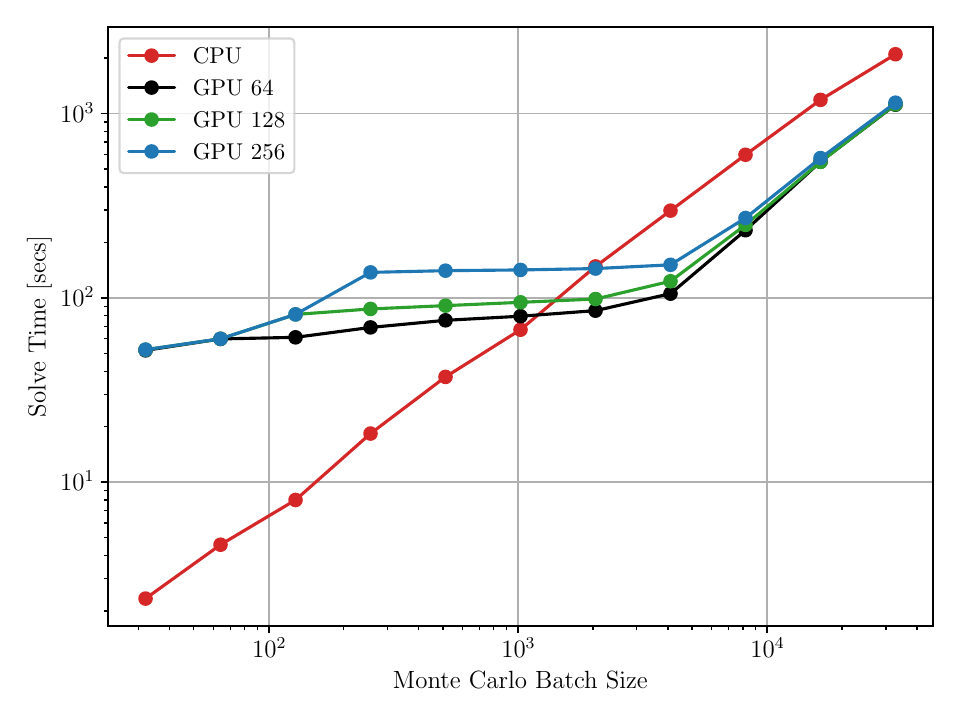}
            \caption{Log-log timing plot}\label{fig:mc-timing-log-log}
        \end{subfigure}
        \begin{subfigure}{.5\textwidth}
            \centering
            \includegraphics[width=\textwidth]{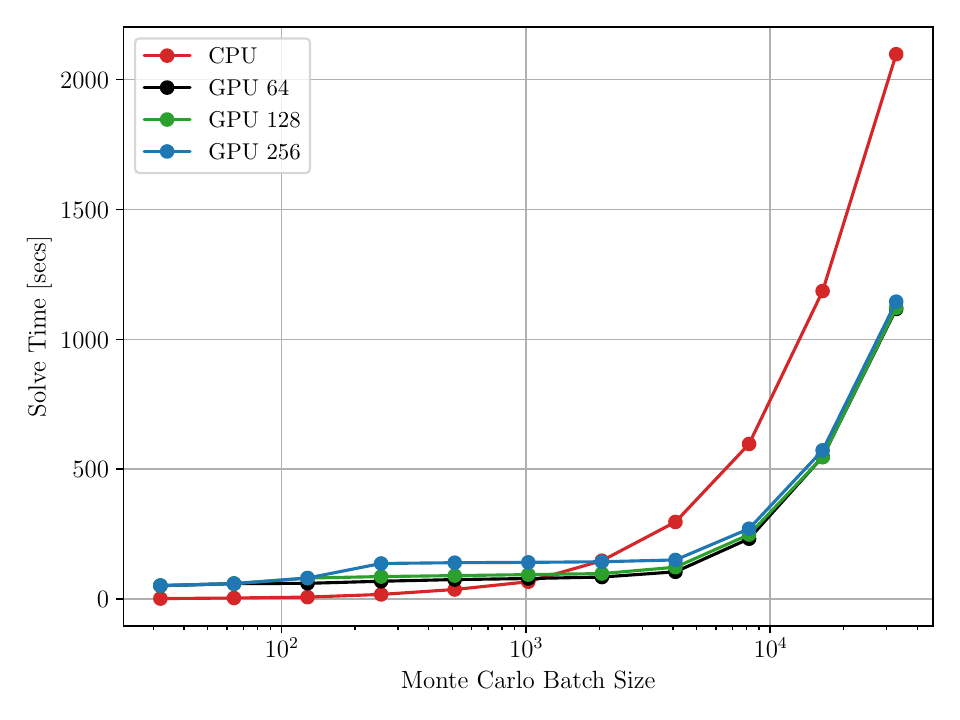}
            \caption{Semi-log timing plot}\label{fig:mc-timing-semi-log}
        \end{subfigure}
        \caption{CPU and GPU Monte Carlo solve-times. GPU results are plotted for the indicated thread-block sizes.}
        \label{fig:mc-timing}
    \end{figure}
Figure \ref{fig:mc-timing} depicts timing results for serial CPU and parallel GPU Monte Carlos for varying Monte Carlo batch sizes and for various thread-block sizes. For the GPU simulations, only the solve-time is reported since execution times for \texttt{malloc} and \texttt{memcpy} were an order of magnitude smaller than the solve-time.

As mentioned in Section \ref{section:gpu-architecture}, calling a kernel requires two parameters: thread-block size and grid size, which is the number of thread-blocks in the grid. The product of these two parameters dictates the number of threads launched. This gives users a choice in choosing the thread-block and grid size when launching a kernel. For example, if we wanted to perform a Monte Carlo with a batch size of 512, we could choose a thread-block size of 128 and grid size of four or a thread-block size of 256 and grid size of two. Here, we perform Monte Carlo runs on the GPU with maximum thread-block sizes of 64, 128, and 256. For batch sizes below the maximum thread-block size, we set the thread-block size to be the batch size, and the grid size to one. For example, if we had a maximum thread-block size of 64 and were conducting a Monte Carlo run with a batch size of four, we would set the thread-block size to 4 and grid size to one. On the other hand, if we were conducting a Monte Carlo run with a batch size of 128, we would set the thread-block size to 64 and the grid size to two.

From Figure \ref{fig:mc-timing}, it is evident that the Monte Carlo run on the GPU is slower than that on the CPU for small batch sizes, and faster for larger batch sizes, with a crossover point roughly between 1024 and 2048 Monte Carlo runs. Since CPUs, unlike GPUs, are optimized for serial tasks such as SCP, it is unsurprising that for smaller runs, the CPU is faster. However, as the batch size increases, the GPU is able to leverage parallelism to provide acceleration. In other words, although the GPU takes longer than the CPU to solve a single SCP problem, it is able to accelerate the Monte Carlo run when given a large number of problems to solve.

The log-log timing plot for the CPU implementation, shown in Figure \ref{fig:mc-timing-log-log}, is linear with a constant slope, as expected, since doubling the number of runs will double the solve-time for a serial implementation. However, the timing plots for the GPU exhibit three phases as the batch size increases. Initially, the solve-time increases until the maximum block size is reached, then the solve-time stays relatively constant, and finally, it increases again. In the first phase, since all threads are in a single thread-block, all the computations for the Monte Carlo are being performed by a single streaming multiprocessor. After the number of runs exceeds the maximum size of the thread-block, a second thread-block is populated, which allows a second streaming multiprocessor to compute these trajectories in parallel with the first thread-block. This results in the second phase, wherein the solve-time curve is flat. The number of runs after which the curve flattens out is dictated by the maximum thread-block size. The third phase, where increasing the number of runs increases the solve-time again, is likely due to the GPU not being able to allocate any more streaming multiprocessors to the thread-blocks, forcing them to execute thread-blocks serially.

%% file: sections/8-conclusion.tex
\FloatBarrier
\section{Conclusion}
We introduce a GPU-accelerated Monte Carlo framework for general nonconvex trajectory optimization problems using sequential convex programming. By comparing a parallel GPU-based implementation of a Monte Carlo framework for the 6-DoF powered-descent guidance problem to a serial CPU-based implementation, we show that for simulations over $\sim$1000 runs, the GPU-based implementation is faster. Future work will investigate parallelized multi-start for nonconvex trajectory optimization, wherein a number of initial guesses are provided, the resulting problems are solved in parallel, and the converged trajectory with the lowest cost is chosen.